\documentclass[final,3p, times, 12pt]{elsarticle}

\usepackage{amssymb,amsmath,color}
\usepackage[mathscr]{eucal}
 \usepackage{lineno}

 \journal{ }
\renewcommand{\journal}{  TBD }
\newtheorem{thm}{Theorem}

\newtheorem{lemma}[thm]{Lemma}

\newtheorem{defin}{Definition}
\newproof{pf}{Proof}

\begin{document}

\begin{frontmatter}



\title{An Erd\H{o}s-Ko-Rado theorem in general linear groups}

\author[JG]{Jun Guo}\ead{guojun$_-$lf@163.com}
\author[KW]{Kaishun Wang\corref{cor}}
\ead{wangks@bnu.edu.cn}
\cortext[cor]{Corresponding author}

\address[JG]{Math. and Inf. College, Langfang Teachers'
College, Langfang  065000,  China }
\address[KW]{Sch. Math. Sci. \& Lab. Math. Com. Sys.,
Beijing Normal University, Beijing  100875, China}

\begin{abstract}
Let $S_n$ be the symmetric group on $n$ points. Deza and Frankl
[M. Deza and P. Frankl, On the maximum number of permutations with
given maximal or minimal distance, J. Combin. Theory Ser. A  22
(1977) 352--360] proved that if ${\cal F}$ is an  intersecting set
in $S_n$ then $|{\cal F}|\leq(n-1)!$. In this paper we consider
the $q$-analogue version of this result. Let $\mathbb{F}_q^n$ be
the $n$-dimensional row vector space over a finite field
$\mathbb{F}_q$ and $GL_n(\mathbb{F}_q)$ the general linear group
of degree $n$. A set ${\cal F}_q\subseteq GL_n(\mathbb{F}_q)$ is
{\it intersecting} if for any $T,S\in{\cal F}_q$  there exists a
non-zero vector $\alpha\in \mathbb{F}_q^n$ such that $\alpha
T=\alpha S$. Let ${\cal F}_q$ be an  intersecting set in
$GL_n(\mathbb{F}_q)$. We show that $|{\cal F}_q|\leq
q^{(n-1)n/2}\prod_{i=1}^{n-1}(q^i-1)$.
\end{abstract}

\begin{keyword}
Erd\H{o}s-Ko-Rado theorem\sep general linear group



\end{keyword}
\end{frontmatter}

\section*{}
The Erd\H{o}s-Ko-Rado theorem \cite{EKR} is a central result in extremal combinatorics.
There are many interesting proofs and extensions of this theorem,
for a summary see \cite{Deza2}.

Let $S_n$ be the symmetric group on $n$ points. A set ${\cal
F}\subseteq S_n$ is {\it intersecting} if for any $f, g\in{\cal
F}$ there exists an $x\in [n]$ such that $f(x) = g(x)$. The
following result is an Erd\H{o}s-Ko-Rado theorem for intersecting
families of permutations.

\begin{thm}\label{thm1.3}
Let ${\cal F}$ be an  intersecting set in $S_n$. Then
\begin{itemize}
\item[\rm(i)]
{\rm(Deza and Frankl \cite{Deza})} $|{\cal F}|\leq(n-1)!$.

\item[\rm(ii)]
{\rm(Cameron and Ku \cite{Cameron})} Equality in {\rm(i)} holds if and only if
${\cal F}$ is a coset of the stabilizer of a point.
\end{itemize}
\end{thm}

 Wang and Zhang \cite{WZ}
gave a simple proof of Theorem~\ref{thm1.3}. Recently, Godsil and
Meagher \cite{Godsil} presented another proof.

In this paper we consider the $q$-analogue of
Theorem~\ref{thm1.3}, and obtain an Erd\H{o}s-Ko-Rado theorem in
general linear groups.

Let $\mathbb{F}_q$ be a finite field and $\mathbb{F}_q^n$ the
$n$-dimensional row vector space over $\mathbb{F}_q$. The set of
all $n\times n$ nonsingular matrices over $\mathbb{F}_q$ forms a
group under matrix multiplication, called the {\it general linear
group} of degree $n$ over $\mathbb{F}_q$,   denoted by
$GL_{n}(\mathbb{F}_q)$.  There is an action of
$GL_{n}(\mathbb{F}_q)$ on $\mathbb{F}_q^{n}$ defined as follows:
\begin{eqnarray}
\mathbb{F}_q^{n}\times GL_{n}(\mathbb{F}_q)&
\longrightarrow&
\mathbb{F}_q^{n}\nonumber\\
((x_1,x_2,\ldots,x_{n}),T)&\longmapsto&
(x_1,x_2,\ldots,x_{n})T.\nonumber
\end{eqnarray}
Let $P$ be an $m$-subspace of $\mathbb{F}_q^{n}$.
Denote also by $P$ an $m\times n$ matrix of rank $m$ whose rows
span the subspace $P$ and call the matrix $P$ a matrix
representation of the subspace $P$.

\begin{defin}
A set ${\cal F}_q\subseteq GL_n(\mathbb{F}_q)$ is {\it intersecting}
if for any $T, S\in{\cal F}_q$ there exists a non-zero vector
$\alpha\in \mathbb{F}_q^{n}$ such that $\alpha T = \alpha S$.
\end{defin}

In this paper, we shall prove the following result:

\begin{thm}\label{thm1.4}
Let ${\cal F}_q$ be an  intersecting set in
$GL_n(\mathbb{F}_q)$. Then $|{\cal F}_q|\leq
q^{(n-1)n/2}\prod_{i=1}^{n-1}(q^i-1)$.
\end{thm}

For the group $GL_n(\mathbb{F}_q)$ we can define a graph, denoted
by $\Gamma$, on vertex set $GL_n(\mathbb{F}_q)$ by joining $T$ and
$S$ if they are intersecting. Since $GL_n(\mathbb{F}_q)$ is an
automorphism group of $\Gamma$, this graph is vertex-transitive.

In order to prove Theorem~\ref{thm1.4}, we require a useful lemma
obtained by Cameron and Ku and a classical result about finite geometry.

\begin{lemma}{\rm (\cite{Cameron})}\label{lem2.1}
Let $C$ be a clique and $A$ a coclique in a vertex-transitive
graph on $v$ vertices. Then $|C||A|\leq v$. Equality implies that
$|C\cap A|=1$.
\end{lemma}

An {\it $n$-spread} of $\mathbb{F}_q^l$ is collection of
$n$-subspaces $\{W_1,\ldots,W_t\}$ such that every non-zero vector
in $\mathbb{F}_q^l$ belongs to exactly one $W_i$.

\begin{thm} {\rm(\cite{Demb})}\label{Thm:spread}
An $n$-spread of $\mathbb{F}_q^l$ exists if and only if $n$ is a
divisor of $l$.
\end{thm}

\begin{lemma}\label{lem2.4}
Let $\alpha(\Gamma)$ be the size of the largest coclique of
$\Gamma$. Then $\alpha(\Gamma)=q^n-1$.
\end{lemma}

\begin{pf} By Theorem~\ref{Thm:spread}, there exists an $n$-spread $\{W_0,W_1,\ldots,W_{q^n}\}$ of
$\mathbb{F}_q^{2n}$.  Since $W_0\cap W_{q^n}=\{0\}$ and
$W_0+W_{q^n}=\mathbb{F}_q^{2n}$, by \cite[Theorem~1.3]{wanbook},
there exists a  $G\in GL_{2n}(\mathbb{F}_q)$ such that
$W_0G=(I^{(n)}\;0^{(n)}),W_{q^n}G=(0^{(n)}\;I^{(n)})$, and
$\{W_0G,W_1G,\ldots,W_{q^n}G\}$ is an $n$-spread of
$\mathbb{F}_q^{2n}$, where $I^{(n)}$ is the identity matrix of
order $n$ and $0^{(n)}$ is the zero matrix of order $n$. Without
loss of generality,  we may assume that $W_0=(I^{(n)}\;0^{(n)})$
and $W_{q^n}=(0^{(n)}\;I^{(n)})$. Then each $W_i\,(1\leq i\leq
q^n-1)$ has the matrix representation of the form
$(I^{(n)}\;T_i)$, where $T_i\in GL_n(\mathbb{F}_q)$. For all
$1\leq i\not=j\leq q^n-1$, since $W_i+W_j$ is of dimension $2n$,
$T_i-T_j\in GL_n(\mathbb{F}_q)$.   By the fact that $T_i-T_j\in
GL_n(\mathbb{F}_q)$ if and only if $\alpha T_i\not=\alpha T_j$ for
all $\alpha\in\mathbb{F}_q^n\backslash\{0\}$,
$\{T_1,\ldots,T_{q^n-1}\}$ is a coclique of $\Gamma$; and so
$\alpha(\Gamma)\geq q^n-1$.

Suppose $\alpha(\Gamma)>q^n-1$ and ${\cal
I}=\{T_1,T_2,\ldots,T_{\alpha(\Gamma)}\}$ is a  coclique of
$\Gamma$. Then $T_i-T_j\in GL_n(\mathbb{F}_q)$ for all $1\leq
i\not=j\leq\alpha(\Gamma)$. Take $W_0=(I^{(n)}\;0^{(n)})$,
$W_{\alpha(\Gamma)+1}=(0^{(n)}\;I^{(n)})$ and
$W_i=(I^{(n)}\;T_i)\;(1\leq i\leq \alpha(\Gamma))$. Then $W_k\cap
W_l=\{0\}$ for all $0\leq k\not=l\leq\alpha(\Gamma)+1$. The number
of non-zero vectors in
$\bigcup_{k=0}^{\alpha(\Gamma)+1}W_k\subseteq \mathbb{F}_q^{2n}$
is
 $(\alpha(\Gamma)+2)(q^n-1)>(q^n+1)(q^n-1)=q^{2n}-1$, a
contradiction. \qed \end{pf}

Combining  Lemma~\ref{lem2.1} and Lemma~\ref{lem2.4}, we complete
the proof of Theorem~\ref{thm1.4}.

Let $G_v$ be the stabilizer of a given non-zero vector $v$ in
$GL_n(\mathbb F_q)$.    Then $G_v$ is an intersecting set meeting
the bound in Theorem~\ref{thm1.4}. It seems to be interesting to
characterize the intersecting sets meeting the bound in
Theorem~\ref{thm1.4}.

\section*{Acknowledgment}

This research is partially supported by    NSF of China (10971052,
10871027),   NCET-08-0052, Langfang Teachers' College (LSZB201005),
and   the Fundamental Research Funds for
the Central Universities of China.

\end{document}